\documentstyle[12pt]{article}
\newcounter{const}

\begin{document}
%\begin{frontmatter}
%\title{A Remark on the  Inequalities of Bernstein - Markov Type in Exponential 
% Orlicz and Lorentz Spaces.}
%\author[Ostrovsky]{E.Ostrovsky}
%\address[Ostrovsky]{Department of Mathematic and Computer Science, 84105, Ben - Gurion 
%University of the Negev, Beer - Sheva, Ben Gurion street, 2; P.O. BOX 61; Israel; \\
% e - mail: galaostr@cs.bgu.ac.il; \ \ galo@list.ru}

%\begin{abstract} 
% We prove in this article the generalizations on the exponential Orlicz and Lorentz
% spaces  inequalities of Bernstein - Markov type for algebraic polynomials and rational 
% functions.
%\end{abstract}

%\begin{keyword}
% Polynomials, Exponential Orlicz Spaces, Markov's and Bernstein's Inequalities, 
% Equivalent Norms.
%\end{keyword}

%\end{frontmatter}

 \large

\begin{center}

{\sc A REMARK ON THE INEQUALITIES OF BERNSTEIN - MARKOV TYPE 
IN EXPONENTIAL ORLICZ AND LORENTZ SPACES.}\\
\vspace{1mm}
{\bf E.I.OSTROVSKY.} \footnote{ Department of Mathematics and Computer science, 
Ben - Gurion University, Beer - Sheva, Israel. E - mail: galaostr@cs.bgu.ac.il} \\
\vspace{1mm}

\end{center}

{\bf Abstract.} {\it We prove in this article the generalizations on the exponential 
Orlicz spaces Markov's - Bernstein's inequalities for  algebraic
polynomials and rational functions.}\\

\vspace{2mm}


\vspace{2mm}


\section {Introduction. Statement of problem.}\par
 
 It is well - known the application to the approximation theory Markov's (or Bernstein's) 
inequalities, for example, for inverse theorems of approximation theory; 
 see for instance ( [7], p. 208;  [13], [14])  etc. \par

 Let us denote by $ {\bf A} $ 
the set of all algebraical polynomials defined on the $ x \in X, \ X = [-1,1]; 
$ by $ {\bf T } $ the set of all 
trigonometrical polynomials on the set $ X = [0, 2 \pi], $ and by $ {\bf R} $ the set of all 
rational functions defined on the  $ X = [-1,1] $ without poles on $ X; $ for $ Q \in {\bf A} 
\cup {\bf T} $ the symbol $ \deg Q $ will denote the usually degree of $ Q; $ for 
irreducible fraction 
$ Q(x) = Q_1(x)/Q_2(x), \ Q_2(x) \ne 0, \ x \in [-1,1], Q_{1,2} \in {\bf A} $ we define 
$ \deg Q = \max(\deg Q_1, \deg Q_2). $  The space $ X $ is equipped usually {\it normed }
Lebesgue measure $ \mu(dx)= Cdx; \ \mu(X) = 1. $  \par

{\bf T}. For all $ Q \in {\bf T} $ 
 hold the famous generalized Bernstein's - Zygmund's inequalities: for all rearrangement 
invariant space $ S $ on the set $ X $
$$
||dQ/dx||S \le \deg Q \cdot ||Q||S, \eqno(1.1)
$$
([12], p. 36;  [20])  and for $ p \in (0,1) $
$$
||dQ/dx||_p \le \deg Q \cdot ||Q||_p,
$$
([7], p. 104),  where by definition for all values $ p > 0 $

$$
||f||_p \ = \ I^{1/p}(|f|^p), \ \ I(f) \stackrel{def}{=} \int_X f(x) \ \mu(dx).
$$

{\bf A}. For all $ Q \in {\bf A}, \ p > 0 $

$$
||dQ/dx||_p \le C^{1 +1/p} \cdot (\deg Q)^2 \cdot ||Q||_p,  
$$
where the symbol $ C $ will denote  (here and further) 
 some  {\it absolute } constant ([4], p.406). In the case $ p > 2 $ 

$$
||dQ/dx||_p \le K(p) \ (\deg Q)^2 \cdot || Q||_p, 
$$

$$
K(p) = \left[ 4 \pi (p+3)^2 /(p \sin(2 \pi/p)) \right]^{1/p}. \eqno(1.2)
$$
 (see [1]).
 Note than  $ p \ge 4 \ \Rightarrow K(p) \le K(4) = (49 \pi)^{1/4} \approx 3.52238228\ldots \ . $ \par
{\bf R}. Let now $ r = 1,2,\ldots, \ p \in (0,\infty), \gamma = \gamma(p,r) =
p/(pr+1), \ Q \in {\bf R}. $ Then
$$
||Q^{(r)}||_{\gamma} \le D_1(p, r) \cdot (\deg Q)^r \cdot ||Q||_p, \eqno(1.3)
$$
where at $ \ p \ge 4 \ \Rightarrow $
$$
D_1(p,r) \le D(r) \stackrel{def}{=} \exp(1/e) \cdot r! \cdot (4/3)^{r+1/4} \cdot (r+1/4)^{r+1/4}, 
$$
see [13], [14]; in [12], p. 300 - 301 was proved that this estimation (1.3) is exact in 
different senses. \par
 For the nonnegative General Algebraic Polynomials (GAP), i.e. for the functions of a view 
$$
Q(x) = \prod_{j=1}^m |x-z(j)|^{r(i)}, \ r(j) \ge 1, \  z(j) = a(j) + b(j)i, \ i = \sqrt{-1},
$$
$$
 \deg Q \stackrel{def}{=} \sum_{j=1}^m r(j), \ x \in [-1,1],
$$
in [4], p.402 - 406 was proved the inequality:

$$
||dQ/dx||_p \le C^{1+1/p}  \cdot (\deg Q)^2 \cdot ||Q||_p, \ p > 0. \eqno(1.4)
$$

 There are many generalizations of this inequalities on the Muntz polynomials [10], spline 
functions [15] and so one ([2], [3], [9], [10], [11], [19]). \par
{\bf Our goal is generalization of inequalities (1.2) and (1.3) on the exponential Orlicz 
spaces, i.e. when in the left and right sides of inequalities (1.2), (1.3) 
instead $ L_p $ norms are the norms on some Orlicz spaces. } \par

\section { Description of using Orlicz spaces.}\par
 We will consider here a so - called exponential Orlicz spaces on the  $ X $ with usually 
 {\it normed} Lebesgue measure $ \mu. $  Recall here  that  if the function 
$ N = N(u),  u \in R^1  $ is some $ N \ - $ Orlicz function (even, continuous, 
downwards convex, $ N(u) \ge 0, \ N(u) = 0 \ \Leftrightarrow u = 0, $
strictly increasing in the self - line  $ R^1_+ $), then the Orlicz norm $ ||f||L(N) $
in the Orlicz's space $ L(N) $
of a (measurable) function $ f: R \to R $ relative to the $ N \ - $ Orlicz function $ N = N(u) $
may be defined by the formula
$$
||f||L(N) = \inf \{l, \ l > 0, \ I(|f|/l) \le 1 \}.
$$
As a particular case, if $ N(u) = |u|^p, \ p = const \ge 1 $
we obtain the classical $ L_p = L_p(R) \ $ spaces with norm $ ||f||_p. $ More information about 
Orlicz spaces see in the books [16], [17]. \par
 A very important class of $ N \ - $ Orlicz functions are the so - called 
$ EOF \ = $ Exponential Orlicz Functions  and correspondent Exponential Orlicz Spaces $ EOS $ 
(like to the terminology of articles [6], [8]) will be considered.  We give here a more general
definition of this spaces.\par
 Let $ \varphi = \varphi(z), \ z \ge 0 $ be some continuous  function such that $ \varphi(z) = 0 \
\Leftrightarrow z = 0, $ and the function $ h(y) := \varphi(\exp y), 
\ y \in [-\infty, \infty) $ is strictly increasing, downward 
convex and 
$$
 \sum_{k \ge 3} \exp(h(k) - h(k+1)) < \infty. 
$$
 The set of all those function we will denote $ \Phi; \ \Phi = \{\varphi \}. $
Let us define the following $ N \ - $ Orlicz function $ N = N(u) = N(\varphi;u): $ at $ |u| \le C_1 $
$$
\ N(\varphi;u) = C_2 \ |u|, \ C_1, C_2 = C_{1,2}(\varphi(\cdot)) \in (0,\infty)  \eqno(2.1)
$$
and for $ |u| > C_1 $

$$
 N(\varphi; u) = \exp \varphi(|u|). \eqno(2.2)
$$
and will denote the correspondent Orlicz's space and norm $ B(\varphi) \stackrel{def}{=} 
L(N(\varphi, \cdot)), \\ ||\cdot||B(\varphi): $  \par
$$
||f||B(\varphi) = ||f||L(N(\varphi; \cdot)).
$$
 It is very simple to prove the existence for all $ \varphi \in \Phi $
the values $ C_1= C_1(\varphi), C_2= C_2(\varphi) $ such that $ N(\varphi;u) $ is some $ N \ - $ 
Orlicz's function.\par
 We denote also $ C_3 = C_3(\varphi) = \max(1,C_1, 1/C_2), $ 

$$
k_0=\max \{ 4+ \max(\log C_1,1), h_+^{*/}(1) \},
$$

$$
C_4(\varphi) = e^2 \left[ N(\varphi; \exp(k_0 -2)) + \sum_{k \ge k_0} \exp(h(k-1)- h(k)) \right].
$$

 Let us consider some  examples. Put $ \varphi(z) = \varphi_{m,r}(z) = z^m \  
\log^{-mr} [(\exp(m + |r|) + z], $
$ z \ge 0, \ m = const \in (0, \infty); $  
or $ \varphi(z) = \varphi_{\nu}(z) = \log^{1+ \nu}( 1 + z), z \ge 0, \ \nu = const > 0.  $ 
Then $  \varphi_{m,r}(\cdot) \in \Phi, \ \varphi_{\nu} (\cdot) \in \Phi. $ We will denote the 
norm in this spaces as
$$
 ||f||B(m,r) = ||f||B(\varphi_{m,r}(\cdot)), \ \ ||f||B(m) = ||f||B(m,0).
$$

\section{Main results.}

 We will prove that the direct generalization of inequality (1.2)
is true for exponential Orlicz spaces.  Denote for $ \varphi \in \Phi $
$$
W(\varphi;n) = \sup_{Q \in A, Q \ne 0} \frac{||dQ/dx||B(\varphi)}{||Q||B(\varphi)}, \eqno(3.1)
$$
where  $ "\sup" $ is calculated over all the algebraic polynomials $ Q; \ Q \ne 0, 
\deg Q = n. $ \\ 

{\bf Theorem 1.} {\it For all} $ \ \varphi(\cdot) \in \Phi $ {\it there exists } $ C_5(\varphi) 
\in (0,\infty) $ such that 
$$
C_5(\varphi) \ n^2 \le W(\varphi,n) \le n^2 \cdot K(4) \cdot \max(1, \psi(4))  \cdot C_{4} \cdot C_3.
\eqno(3.2)
$$
 
 Note than since $ \mu(X) = 1, $ the inequalities (3.2) hold for all the Orlicz's spaces with the 
$ N \ - $ functions which are equivalent to $ N(\varphi; u). $ \par
  In order to formulate an other result, we introduce some new notations.
 For  $ \varphi \in \Phi $  we denote
$$
\psi(p) = \psi(\varphi; p) = \exp(h^*(p)/p),
$$
 where  $ h^*(p) $  denotes the classical Young - Fenchel, or Legendre transform:

$$
h^*(p) = \sup_{y \in (-\infty, \infty) } (py - h(y)),
$$
 We  define for  $ r = 1,2,\ldots, \ \varphi \in \Phi $ a new quasinorm

$$
||f||V(\varphi;r) \stackrel{def}{=} \sup_{\beta \in (4/(4r+1), 1/r)} ||f||_{\beta}/
\psi(\beta/(1-r\beta))
$$
and the correspondent space of measurable  functions $ V(\varphi; r) $ with finite norm 
$ ||f||V(\varphi;r) < \infty. $ \\

{\bf Theorem 2. } $ \forall \varphi \in \Phi, \ \forall Q \in {\bf R} $

$$
 ||Q^{(r)}||V(\varphi;r) \ \le \ C_4(\varphi) \cdot D(r) \cdot (\deg Q)^r \cdot
||Q||B(\varphi). \eqno(3.3)
$$

\section{ Auxiliary result.}
  
 Let us introduce a {\it new } Banach space $ G( \varphi), \ \varphi \in \Phi, $ 
as a set of all measurable functions $ f:X \to R  $ with finite norm

$$
||f||G( \varphi) \stackrel{def}{=} \sup_{p \ge 1} |f|_p /\psi(p) < \infty.
$$

 Note than by virtue of Iensen - Lyapunov inequality 
$$
\sup_{p \ge 4} ||f||_p/\psi(p) \le ||f||G(\varphi) \le \max(1, \psi(4)) \ 
\sup_{p \ge 4} ||f||_p /\psi(p).
$$

{\bf Theorem 3.} 
 {\it We propose that the norms } $ || \cdot ||B( \varphi) $ {\it and } 
 $ || \cdot ||G(\varphi) $ {\it are equivalent: } 

$$
C_3^{-1} \ ||f||G( \varphi) \le ||f||B( \varphi) \le C_4 \ ||f||G( \varphi). \eqno(4.1)
$$

{\bf Proof of theorem 3.} 
 Assume at first $  ||f||B( \varphi) < \infty. $ Without loss of generality we 
can suppose  

$$
I (N(\varphi; |f|) = 1. 
$$
 Let us introduce the function

$$
g(p) = \sup_{z > 0} z^p / N( \varphi; z).
$$
 We have for all the values $ p \ge 1: $ 
$$
g(p) \le  \max \left[\max_{z \in (0, C_1] } C_2^{-1}  z^{p-1},  \ \sup_{z \ge C_1} z^p \ 
\exp (- \varphi(z)) \right]  \le
$$

$$
\max \left[ C_2^{-1} \max(C_1,1)^{p-1}, \ \exp( \sup_{z \ge C_1} (p \log z - \varphi(z))) \right] <
$$
$$
 \max \left[ C_2^{-1} \max(C_1,1)^{p-1}, \ \exp(\sup_{v \in (-\infty,\infty)} 
(p v - h(y))) \right] \le
$$

$$
 \max \left[C_3^{p}, \  \exp h^*(p) \right] \le C_3^p \exp h^*(p).
$$

 Following, for the values $ p \ge 1 $ we have: $ z \ge 0 \ \Rightarrow $

$$
z^p \le g(p) \ N(\varphi; z) \le C_3^p \ \psi^p(p) \ N(\varphi; z). 
$$

Therefore
$$
 |f|^p \le C_3^p \ \psi^p(p) \ N(\varphi; h^*(|f|)), \ \ 
|f|_p \le C_3( \varphi) \ \psi(p),
$$

$$
||f||G( \varphi) \le C_3( \varphi(\cdot)) < \infty.
$$

 Inverse, suppose 

$$
|f|_p^p  \le \exp \left(h^*(p) \right), \ p \ge 1.
$$
 We have by virtue of Chebyshev's inequality for the values $ w \ge e^2:$

$$
T(|f|,w) \stackrel{def}{=} \mu \ \{x: |f(x)| > w \} \le \exp \left(h^*(p) - p \log w \right),
$$
and after the minimization over $ p: $
$$
T(|f|,w) \le \exp \left( -\sup_{p \ge 1} (p \log w - h^*(p) \right) \le 
$$
$$
\exp \left(- \sup_p (p \log w - h^*(p) \right) =
$$ 

$$
\exp \left(-h^{**}(\log w) \right) = \exp \left( -h( \log w) \right)
$$
by virtue of theorem Fenchel - Moraux. \par
 We conclude for 
the value  $ \varepsilon = \exp(-2), $ choosing $ W(k) = \exp(k) $ and denoting

$$
U(k) = U(|f|,k) = \{x: \ W(k) \le |f(x)| < W(k+1) \}: 
$$

$$
I(\exp(N(\varphi; \varepsilon |f|))) =\int_{ \{x: |f(x)| \le \exp(k_0) \} } N(\varepsilon |f(x)|)
 \ d \mu \ +
$$

$$
\int_{ \{x: |f(x)| > \exp(k_0) \} } N(\varphi;\varepsilon |f(x)|) \ d \ \mu \le 
N(\varphi;\exp(k_0-2)) \ +
$$

$$ 
\sum_{k \ge k_0} \int_{U(k)}\exp(\varphi(\varepsilon |f(x)|) ) \ dx \le 
$$

$$
N(\varphi; \exp(k_0-2) +\sum_{k \ge k_0} \exp 
\left[\left( h(\varepsilon \  W(k+1) \right) \right] \ \cdot \ \left[ T(|f|, W(k)) \right] \le
$$

$$
N(\varphi; \exp(k_0-2)) + \sum_{k \ge k_0} \exp(h(k) - h(k+1)) =C_4 \ e^{-2} < \infty.
$$
This completes the proof of theorem 3.\par
 Note than this result is some generalization of [5], p. 309 - 314, [6], [8], [17], p. 305.\par
 For example, let $ N(u) = N_{m,r}(u)= N(\varphi_{m,r}(\cdot);u). \  $  
It follows from theorem 3 that 
$$
 ||f||L(B(m,r))  < \infty \ \Longleftrightarrow \ \sup_{p \ge 4} |f|_p \ \cdot 
 \left(p^{-1/m} \log^{m r} p \right) < \infty,
$$
 or equally

$$
\exists \varepsilon > 0, \ I(N(\varphi_{m,r}; \varepsilon |f|) < \infty \ \Longleftrightarrow
\ \sup_{p \ge 4} |f|_p \ \cdot \left( p^{-1/m} \ \log^{mr} p \right) < \infty.
$$
{\bf Notice.} Let us introduce the {\it weight Lorentz } norm:

$$
||f||^*_b G(\varphi) = \sup_{p \ge 1} ||f||_{p,b}/\psi(p),
$$
where $ ||f||_{p,b} $ is the Lorentz norm (more exactly, seminorm):

$$
||f||_{p,b} = \left[\int_0^{\infty} T^{p/b}(|f|,x) \ dx^b \right]^{1/b}, 
$$
$ p \in [1, \infty), \ b \in [1,\infty], $ where in the case $ b = +\infty $ 
$$ 
||f||_{p,\infty} = \sup_{x \ge 0} \left( x \ T^{1/p}(|f|,x) \right).
$$
 Using the embedding theorem for Lorentz spaces it is easy to prove as well as by 
proving of theorem 4 that all the following norms are equivalent
$$
||\cdot||B( \varphi) \ \sim \ ||\cdot|| G( \varphi) \ \sim \ ||\cdot||^*_b G(\varphi)
$$ 
 with constants does not depending on $ b. $ Therefore, it is easy to formulate theorems 
1,2 in the terms of those spaces.\par

\section{ Proof of the main results.}
{\bf Proof of theorem 1.} The low bound in (3.2) is attained, for instance, on the 
so - called Jacobi ultraspherical polynomials $ Q = P_n^{2,2}: $
$$
P_n^{2,2}(x) = \frac{(-1)^n (1-x^2)^{-2}}{2^n \ n!} \left( \frac{d}{dx} \right)^n 
\left[(1-x^2)^{n+2} \right]:
$$

$$
||(d/dx) \ P_n^{2,2}||_p \ge C \ n^2 \ ||P_n^{2,2}||_p, \ n = \deg P_n^{2,2}. \eqno(5.1)
$$

see [18], p. 66, 165;  [11]. The upper bound may be simple prove by virtue of theorem 3. 
Namely, suppose $ Q \in A, \ Q \ne 0, \ \deg Q = n \ge 2, \ ||Q||B(\varphi) = 1.$ Then 

$$
||Q||G(\varphi) \le C_3||Q||B(\varphi) = C_3 < \infty;
$$

$$
 \sup_{p \ge 4} \ \left( ||Q||_p / \psi(p) \right) \le C_3; \ \Rightarrow  ||Q||_p \le 
C_3 \psi(p) 
$$
for all the values $ p \ge 4.$ It follows from [1] that
$$
||dQ/dx||_p \le n^2 \cdot K(4) \cdot \psi(p) \cdot C_3. 
$$
 Therefore 
$$
||dQ/dx||G(\psi) \le n^2 \ K(4) \ \max(1, \psi(4)) \cdot C_3, 
$$
and finally

$$
||dQ/dx||B(\varphi) \le n^2 \cdot K(4) \cdot C_3 \cdot C_4 \cdot \max(1, \psi(4)) \cdot 
||Q||B(\varphi).
$$

{\bf Proof of theorem 2.} Let $ Q \in R, \ \deg Q = n, \ ||Q||B(\varphi) = 1. $ By virtue of 
theorem 3 we receive:

$$
||Q||G(\varphi) \le C_4; \ \Rightarrow ||Q||_p \le C_4 \ \psi(p), \ p \ge 4.
$$
 From (1.3) it follows
$$
||Q^{(r)}||_{\gamma} \le n^r \ D(r) \ C_4(\varphi) \ \psi(p).
$$
 Since $ p \ge 4, \ \gamma \in [4/(4r+1),1/r). $ After the substitution 
$ \beta = p/(pr+1) \in [4/(4r+1), 1/r), \ $ or $ p= \beta/(1-\beta \ r), $ we obtain
that for all $ \beta \in [4/(4r+1), 1/r) $
$$
||Q^{(r)}||_{\beta} \le n^r \ C_4(\varphi) \ D(r) \ \psi(\beta/(1-\beta \ r)),
$$
which is equivalent to the statement (3.3) of theorem 2. \par 

\section{Concluding Remark.} 
The spaces $ V(\varphi, r) $ are only Frech\'et spaces. Now we 
investigate  the connection between the norm $ || f ||V(\varphi_{m,0};r) $ 
in those spaces ant tail behavior of distribution $ T(|f|,u). $ At first assume that 
$ ||f||V(\varphi_{m,0},r) = 1 $  for some $ m > 0. $
then for all values $ \beta \in [4/(4r+1), 1/r) $
$$
||f||_{\beta} \le \left( \frac{\beta}{1-r \beta} \right)^{1/m}, 
\ \ I \left(|f|^{\beta} \right) \le \left( \frac{\beta}{1-r\beta}  \right)^{\beta/m}.
$$
We obtain using the Chebyshev's inequality for sufficiently large values $ u 
\ge 3: $

$$
T(|f|, u) \le u^{-\beta} \ \left( \frac{\beta}{1-r \beta} \right)^{\beta/m}.
$$
We conclude after the minimization over $ \beta: $
$$
T(|f|,u) \le C_9(m,r) \ u^{-1/r} \ \left[\log u \right]^{1/(mr)}, \ u \ge 3. \eqno(6.1)
$$
 Inverse, if inequality (6.1) holds, then 

$$
I \left( |f|^{\beta} \right) \le C_{11}(\beta,r,m) \int_3^{\infty} x^{\beta-1-1/r} \ 
(\log x)^{1/mr} \ dx \le 
$$

$$
C_{12} (1/r - \beta)^{-(mr+1)/(mr)}; \ \ ||f||_{\beta} \le C_{13} (1/r - \beta)^{(mr+1)/m};
$$

$$
||f||V(\varphi_{m/(mr+1),0}; \ r) \le C_{14}(m,r) < \infty.
$$
\vspace{2mm}
{\bf Aknowledgements}. I am very grateful to prof. V.Fonf, M.Lin (Ben Gurion University, 
Beer - Sheva, Israel) for useful support of these investigations. \par
 This investigation was partially supported by ISF(Israel Science Foundation), grant 
$ N^o $ 139/03.\par

\newpage

{\bf References.} \\

1. M. Baran. New Approach to Markov Inequality in $ L^p. $ In: Approximation Theory;
  Memory of A.K.Varma, p. 75 - 85.  Edited by N.K.Govie. New - York, Basel. 2004. \\
2. D.Benko, and Tama's Erde'lyi. Markov Inequality for Polynomials of degree $ n $ with $ m $ 
distinct Zeroes. Journal Appr. Theory, 2003, 122,  241 - 248.\\
3. D. Benko, Totok Vilmos. Set with Interior Extremal Points for the Markov Inequality. 
Journal of Appr. Theory, 2001, 112, 171 - 188.\\
4. P. Borwein, Tama's Erde'lyi. Polynomials and Polynomial Inequalities. Springer 
Verlag, New York, Berlin, Heidelberg, 1995. \\
5. Buldygin V.V., Mushtary D.M., Ostrovsky E.I., Puchalsky A.I. New Trends in Probability Theory
and Statistics. 1992, Springer Verlag, New York - Berlin - Heidelberg - Amsterdam.\\
6. Cruz David - Uribe, SFO and Miroslav Krbec. Localization and Extrapolation in Orlicz - 
Lorentz Spaces. In: coll. works "Function Spaces, Interpolation Theory and Related Topics".
Editors: M. Cwikel, M. Englis etc., Berlin, New York; 2002, 273 - 283. \\
7. DeVore Ronald A., Lorentz George G. Constructive Approximation. Springer Verlag,
New York - Berlin - Heidelberg, 1993. \\
8. Edmunds D.E., and Krbec M. On decomposition in Exponential Orlicz Spaces. Math. Nachr.,
2000, 213, 77 - 88.\\
9. E. Hille, G. Szeg\"o, J. Tamarkin. On some generalization of a 
theorem of A.Markoff. Duke Math. J. 1937, 3, 729 - 739.\\
10. Tama's Erde'lyi. Markov Type Inequality for Products of Muntz Polynomials. Journal of 
Appr. Theory, 2001, 112, 171 - 188.\\
11. G.K. Kristiansen. Some inequalities for algebraic and trigonometric polynomials. J. 
London Math. Soc. 1979, 20, 300 - 314.\\
12. Lorentz Georg G., Golitsek Manfred V., Makovoz Yuly. Constructive Approximation.
Springer Verlag, New York - Berlin - Heidelberg, 1996. \\
13. A.A. Pekarskii. Estimates of derivatives of rational function in $ L_p[-1,1]. $
Math. Zam., 1986, 39, 388 - 394. \\
14. A.A. Pekarskii. Inequalities of Bernstein type for derivatives of rational function, 
and inverse theorem of rational approximation. Math. Sb., 1984, 124, 571 - 588.\\
15. Pincus A., Schoenberg I.J. Generalized Markov - Bernstein Type Inequalities for Spline 
Functions. In: Studies in spline Functions and Approximation Theory. Academic Press Inc.
 Editors S. Karlin, G. Micchelli, 1976. \\
16. M.M.Rao, Z.D.Ren. Theory of Orlicz Spaces. Marcel Deccer inc., New York, Basel, 1993.\\
17. M.M.Rao, Z.D.Ren.  Applications of Orlicz Spaces.  Marcel Dekker Inc., New York,
Basel, 2002.\\
18. Szeg\"o G. Orthogonal Polynomials. Orthogonal Polynomials. Amer. Math. Soc., 
New York, 1939.\\
19. Varma A.K. On Some Extremal Properties of Algebraic Polynomials. J. Appr. Theory, 1992, 69, 
48 - 54.\\
20. Zygmund A. A remark on conjugate functions. Proceedings of the London Math. Soc. 1932, 34,
392 - 400.\\

%\newpage

%\begin{center}

%{\bf  APPLICATION FORM FOR PARTICIPANTS }\\
%\end{center}

%{\bf First name:} Eugene \\

%{\bf Last (family) name:} Ostrovsky \\

%{\bf Institution:} University Ben Gurion, Beer - Sheva city, Israel, Ben Gurion street, 2;
%building 58, Department of Mathematic and Computer Science, PO BOX 61. \\

%{\bf Location:} Israel, 56209, Rehovot city, Shkolnik street, H.5 Ap. 8; \\
%Beer - Sheva, Ben - Gurion University, Dept. of Mathematic and Computer Science.\\

%{\bf Phone numbers.} Home: (972) - (08) - 9 - 45 - 16 - 13, \  Work: (972) -
%(08) - 647 - 78 - 49. \\

%{\bf E - mail address:} galaostr@cs.bgu.ac.il \\

%{\bf Title of report:} Correct Solvability of Non - Linear Ordinary Differential Equations 
% in Orlicz Spaces.\\

%{\bf References:} ??? publications. \\

\end{document}